\documentclass[12pt,a4paper]{article}
\usepackage{amssymb,amsmath}
\usepackage{amsfonts}
\usepackage{mdwlist}
\numberwithin{equation}{section}
\usepackage{a4}
\newcommand{\para}{\par\vspace{.25cm}}

\newtheorem{theorem}{Theorem}
\newtheorem{lemma}{Lemma}
\newtheorem{cor}{Corollary}

\usepackage{multirow}

\begin{document}
	\baselineskip 18pt
	
		\title{\bf Finite semisimple group algebra of a normally monomial group}
	\author{Shalini Gupta {\footnote{Corresponding author}} \\ {\em \small Department of Mathematics,} \\
		{\em \small Punjabi University, Patiala,  India.}\\{\em
			\small email: shalini@pbi.ac.in} \and  Sugandha Maheshwary {\footnote {Research supported by SERB, India (PDF/2016/000731).}}
		\\ {\em \small Department of Mathematical Sciences,}\\{\em\small Indian Institute of Science Education and Research, Mohali,}\\
		{\em \small Sector 81, Mohali (Punjab)-140306, India.}
		\\{\em \small email: sugandha@iisermohali.ac.in}
	}
	\date{}
	{\maketitle}
	
	\begin{abstract}\noindent
	In this paper, the complete algebraic structure of finite semisimple group algebra of a normally monomial group is described. The main result is illustrated by computing the explicit Wedderburn decomposition of finite semisimple group algebras of various normally monomial groups. The \linebreak automorphism groups of these group algebras are also determined.
\end{abstract}\vspace{.25cm}
	 \vspace{.25cm} 

	\section{ Introduction}Let $\mathbb{F}_{q}$ denote the field containing $q$ elements and let $G$ be a finite group of order relatively prime to $q$, so that 
the group algebra $\mathbb{F}_{q}G$ is semisimple. The knowledge of the algebraic structure of $\mathbb{F}_{q}G$ has applications in coding theory and is useful in describing the automorphism group as well as the unit group of 
$\mathbb{F}_{q}G$. This has attracted the attention of several authors \cite{BGP11,BGP,BGP3,bak1,bak2,Broche,neh1,neh2,OltV}. \par
Broche et al. \cite{Broche} gave  description of semisimple group algbra $\mathbb{F}_{q}G,$ when $G$ is an abelian-by-supersolvable group, by computing its
primitive central idempotents and the corresponding simple components, in terms of subgroups of $G$. Basing on the work in \cite{BKP} and \cite{Broche}, a
more precise description of $\mathbb{F}_{q}G$, where $G$ is a finite metabelian group, has been given by Bakshi et al. in a series of papers 
\cite{BGP11,BGP,BGP3}.

\par The present paper is a contribution to the work in same series. Recall that a group $G$ is called \emph{normally monomial}, if every irreducible character of $G$ is \linebreak\emph{ normally monomial}, i.e., induced from a linear character of a normal subgroup of $G$. 
We provide a complete set of primitive central idempotents and the \linebreak Wedderburn decomposition of $\mathbb{F}_{q}G$, when $G$ is a finite normally monomial group. Since metabelian groups are normally monomial \cite{Basm}, this generalises the main \linebreak result of \cite{BGP3}. It may be remarked that normally monomial groups form a \linebreak substantial class of monomial groups \cite{BM3} and the rational group algebra of this class of groups has been studied \cite{BM1}.  

The main result is given in Section 2, after setting up necessary notation and preliminaries  and is illustrated on a family of metabelian groups. In Section 3, we give applications of the main result on certain $p$-groups, which include a family of non-metabelian but normally monomial groups of order $p^{7}$, $p$ prime, $p\geq 5$. It may be pointed out that for most of these group algebras, the \texttt{GAP} package \texttt{Wedderga} \cite{wedd} practically fails to compute the
Wedderburn decomposition. Further, we provide the explicit structure of $\mathbb{F}_{q}G$, for any group $G$ of order $p^{n}$, $p$ prime, $n<5$. For these groups, the  group $\operatorname{Aut}(\mathbb{F}_{q}G)$ of $\mathbb{F}_{q}$-automorphisms of $\mathbb{F}_{q}G$  has also been computed.

	\section*{Notation}
	Throughout the paper, $G$ denotes a finite group, and $\mathbb{F}_{q}G$ denotes the group \linebreak algebra of $G$ over the field $\mathbb{F}_{q}$ containing $q$ elements, $q$ relatively prime to the order $|G|$ of $G$. The notation used are mostly standard and are listed for the ease of reader.\\
	
	$ \begin{array}{ll}

	H \leq (\unlhd ) G  & H {\rm~is~ a ~subgroup~ (normal~subgroup)~of} ~G\vspace{.1cm}\\

	N_{G}(H) & {\rm the~ normalizer~ of~} H{\rm~in~}G,~ H \leq G\vspace{.1cm}\\
	\left[G:H\right] &  {\rm the~ index~ of~ the~ subgroup}~ H ~{\rm in}~  G,~ H \leq G\vspace{.1cm}\\

	H^g & g^{-1}Hg,~ g \in G,~ H \leq G \vspace{.1cm} \\
	
		\operatorname{core}_{G}(H) & \displaystyle\bigcap_{g\in G}H^{g} ,~{\rm the~ largest~ normal~ subgroup~ of ~}G ~{\rm contained~in}~ H,~ H \leq G\vspace{.1cm}\\
		
			\left[h,g\right] & h^{-1}g^{-1}hg,~ g, h \in G\vspace{.1cm} \\
			
			G' & \{	\left[h,g\right] ~|~ g, h \in G\}\vspace{.1cm} \\
			\end{array}$
			$ \begin{array}{ll}
			|X| & {\rm  the~cardinality ~of ~the ~set}~ X\vspace{.1cm}\\	
		\langle X\rangle & {\rm the~  subgroup~ generated ~by~ the ~subset}~ X ~{\rm of}~ G\vspace{.1cm}\\
		
			\ \psi^{G}& {\rm the~ character }~\psi ~{\rm of ~ a ~subgroup~of }~G,~{\rm induced~ to~ }~ G\vspace{.1cm}\\
	\overline{\mathbb{F}}_{q}& {\rm an ~algebraic ~closure~ of}~ \mathbb{F}_{q}\vspace{.1cm}\\
	
	\operatorname{Irr}(G)& {\rm the~ set~ of~ all~ the ~distinct~
		irreducible ~characters ~of}~ G~{\rm over}~\overline{\mathbb{F}}_{q}\vspace{.1cm}\\

	\operatorname{ker}(\chi)& \{ g \in G \, \mid \, \chi(g) = \chi(1)\},~ \chi \in \operatorname{Irr}(G)\vspace{.1cm}\\
\mathbb{F}_{q}(\chi) & {\rm the~ field~ obtained ~ by ~adjoining ~} \chi (g),~g \in G,~{\rm to} ~ \mathbb{F}_{q} ,~\chi \in \operatorname{Irr}(G)\vspace{.1cm}\\
	\operatorname{Gal}(\mathbb{F}_{q}(\chi)/\mathbb{F}_{q}) & ~{\rm the~Galois ~group ~of }~ \mathbb{F}_{q}(\chi)~ {\rm over}~\mathbb{F}_{q}\vspace{.1cm}\\
	e(\chi) & \frac{\chi(1)}{|G|}{\sum _{g \in G}} \chi(g) g^{-1},~ \chi \in \operatorname{Irr}(G)\vspace{.1cm}\\
	e_{\mathbb{F}_{q}}(\chi)&{ \sum_{ \sigma \,\in\,\operatorname{Gal}(\mathbb{F}_{q}(\chi)/ \mathbb{F}_{q})}} e( \sigma\circ \chi),~{\rm the~ primitive~ central ~idempotent ~of}
	~\mathbb{F}_{q}G \vspace{.1cm}\\ & {~\rm associated}~ {\rm to ~the ~character} ~\chi ,~ \chi \in \operatorname{Irr}(G)\vspace{.1cm}\\
a| b ~(a\nmid b) & a {\rm~ divides ~(does ~not~ divide)~} b,~a,b \geq 1\vspace{.1cm}\\
	\varphi & {\rm Euler's~ phi~ function}\vspace{.1cm}\\
	\operatorname{gcd}(k,n)& {\rm the~ greatest ~common~ divisor~ of~the~integers} ~k~ {\rm and}~ n\vspace{.1cm}\\
	S_{n} & {\rm the ~symmetric~ group~ on}~ n~{\rm  symbols},~n\geq 1\\
	\mathbb{Z}_{n}&{\rm the~ cyclic~ group
		~of~ order}~ n,~n \geq 1\\
	\mathrm{SL}_{n}(\mathbb{F}) & {\rm the ~group~ of}~ n\times n~{\rm~ matrices~ of ~determinant~ 1}, ~{\rm over~ the~ field}~\mathbb{F},~n \geq 1\\
	H_{1}
	\rtimes  H_{2}& {\rm the~ split~ extension~ of~ the~ group}~H_{1}~{\rm by~ the~ group}~H_{2}\\
	
	M_{n}(\mathbb{F})& {\rm  the~ring ~of}~ n \times n~ {\rm matrices~ over~ the~ field}~\mathbb{F},~n \geq 1\vspace{.1cm}\\
	M_{n}(\mathbb{F})^{(s)}& M_{n}(\mathbb{F})\bigoplus M_{n}(\mathbb{F})\bigoplus...\bigoplus M_{n}(\mathbb{F}),~{\rm the~direct ~sum ~of}~ s~{\rm copies},~ s\geq 1\vspace{.1cm}\\
	\end{array}
	$
	
	\section{Normally monomial groups}  
	\subsection{Strong Shoda pairs and idempotents}
	For $H\unlhd K \leq G$, define $$\varepsilon (K, H) :=  \begin{cases} \hat{K}, & K = H\\ \prod (\hat{H}-\hat{L}), &  {\rm otherwise} \end{cases},$$ where $\hat{K} :=\frac{1}{|K|} \sum_{k \in K} k$ and $L$ runs over the normal subgroups of $K$ which are minimal over the normal subgroups of $K$ containing $H$ properly. Set $e(G,K,H)$ to be the sum of distinct $G$-conjugates of $\varepsilon(K,H).$  A {\it strong Shoda pair} \cite{Oli} of $G$ is a pair $(K, H)$ of subgroups of $G$ with the property that \begin{quote}  (i) $H \unlhd K  \unlhd N_{G}(H)$; \\ (ii) $K/H$ is cyclic and a maximal abelian subgroup of $N_{G}(H)/H$;\\ (iii) the distinct $G$-conjugates of $\varepsilon(K, H)$ are mutually orthogonal. \end{quote} Two strong Shoda pairs $(K_{1}, H_{1})$ and $(K_{2}, H_{2})$ are said to be {\it equivalent}, if\linebreak $e(G, K_{
		1}, H_{1})$ = $e(G, K_{2}, H_{2})$ and a complete set of representatives from distinct equivalence classes of strong  Shoda pairs of $G$ is called {\it a complete  irredundant set of strong Shoda  pairs} of $G$.
	\vspace{0.5cm}
			
	We recall the algorithm to compute a complete irredundant set of strong Shoda pairs of a finite normally monomial group $G$, as described in \cite{BM1}.
	
	\para Let $\mathcal{N}$ be the set of all normal subgroups of $G$ and for $N\in \mathcal{N}$,  let $A_{N}$ be a \linebreak normal subgroup of $G$ containing $N$ such that $A_{N}/N$ is an abelian normal subgroup of maximal order in $G/N$. Note that the choice of $A_{N}$ is not unique. However, we need to fix one such $A_{N}$. For a fixed $A_{N}$, set \vspace{.35cm}\\ $\begin{array}{lll}
	\mathcal{D}_{N}: & {\rm the ~set ~of ~all ~subgroups~} D ~{\rm of} ~A_{N}~ {\rm containing ~} N~{\rm  such ~that~}
	\operatorname{core}_{G}(D)=N, \\&  A_{N}/D ~{\rm is ~ cyclic ~  and ~  is ~a ~ maximal~ abelian ~ subgroup~  of~}   N_{G}(D)/D.  \vspace{.2cm}\\

	\mathcal{T}_{N}: &  {\rm a ~set~ of ~representatives~ of~ } \mathcal{D}_{N} {\rm ~  under ~the~ equivalence~ relation~ defined~ by} \\ & {\rm
		conjugacy~ of~ subgroups~in~} G. \vspace{.2cm} \\  \mathcal{S}_{N}: & \{( A_{N},D)~|~ D \in \mathcal{T}_{N}\}.\end{array}$
	
	\para \noindent It thus follows that if $N \in \mathcal{N}$ is such that $G/N$ is abelian, then \begin{equation}\label{e13} \mathcal{S}_{N} = \begin{cases} \{(G, N)\}, & {\rm if}~G/N~{\rm is~ cyclic,}\\ \emptyset, & {\rm otherwise.}
	\end{cases} \end{equation}

	\para \noindent Observe that every pair $(A,D)\in\mathcal{S}(G)$, where \begin{equation}\label{E1} \mathcal{S}(G):=\displaystyle\bigcup_{N \in \mathcal{N}}\mathcal{S}_{N},\end{equation} is a strong Shoda pair of $G$. It has been proved (\cite{BM1}, Corollary 1) that $\mathcal{S}(G)$
	is a complete irredundant set of strong Shoda pairs of $G$, if $G$ is a finite normally monomial group. \\
	
	\noindent {\bf{Remark 1.}}
	A crucial observation in the above algorithm to compute $\mathcal{S}(G)$, for a given finite group $G$, is that the choice of $A_{N}$ is irrelevant. For $N \in \mathcal{N}$, let $A_{N}^{(1)}$ be another normal subgroup of $G$ containing $N$ such that $A_{N}^{(1)}/N$ is an abelian normal subgroup of maximal order in $G/N$ and let $\mathcal{D}_{N}^{(1)}$, $\mathcal{T}_{N}^{(1)}$ and $\mathcal{S}_{N}^{(1)}$ be defined corresponding to $A_{N}^{(1)}$. Then, any pair in $\mathcal{S}_{N}^{(1)}$ is equivalent to a pair in $\mathcal{S}_{N}$ and \linebreak vice-versa. This is because, if $(A_{N}^{(1)},D^{(1)})\in \mathcal{S}_{N}^{(1)}$ and $\psi$ is a complex linear character of $A_{N}^{(1)}$ with kernel $D^{(1)}$, then $\psi^{G}$ is irreducible and hence by (\cite{BM1}, Lemma 1), there exists $(A_{N},D)\in \mathcal{S}_{N}$ such that $e_{\mathbb{Q}}(\psi^{G})$, the primitive central idempotent of the rational group algebra $\mathbb{Q}G$, associated to $\psi^{G}$, is given by $e(G,A_{N},D).$ However, in view of (\cite{
Oli}, 
Theorem 2.1), $e_{\mathbb{Q}}(\
	\psi^{G})=e(G,A_{N}^{(1)},D^{(1)}).$ This gives that $(A_{N},D)$ is equivalent to $(A_{N}^{(1)
	},D^{(1)})$. The reverse conclusion holds similarly.\\
	
	For a strong Shoda pair $(K,H)$ of $G$, let $\mathcal{C}(K/H)$ denote the set of $q$-cyclotomic cosets of $\mbox{Irr}(K/H)$ containing its generators, i.e., if $\chi$ is a generator of $\mbox{Irr}(K/H),$ then an element $C$ of $\mathcal{C}(K/H)$ containing $\chi$ is the set $\{\chi,\,{\chi}^{q},\ldots, 
	{\chi}^{q^{o-1}}\},$ where $n=|K/H|$ and $o=o_{n}(q),$ the order of $q$ modulo $n.$ Suppose that $N_{G}(H)$ acts on $\mathcal{C}(K/H)$ by conjugation, i.e., for $g\in N_{G}(H),\,C \in \mathcal{C}(K/H)$ and $\chi \in C,$
	we define ${C}^{g}=\{{\chi}^{g},\,{\chi^{g}}^{q},\ldots,
	{\chi^{g}}^{q^{o-1}}\},$ where $\chi^{g}(k)=\chi(gkg^{-1}).$ Let $\mathcal{R}(K/H)$ denote the set of distinct orbits of $\mathcal{C}(K/H)$  and let
	$E_{G}(K/H)$ be the stabilizer of any $C\in \mathcal{C}(K/H)$ under the above action. It is easy to see that 
	$|\mathcal{R}(K/H)|=\frac{\phi(n)|E_{G}(K/H)|}{|N_{G}(H)|o_{n}(q)}.$ 
	
	\par For $C\in \mathcal{C}(K/H)$ and $\chi\in C,$ set 
	$\varepsilon_{C}(K, H) =\frac{1}{|K|}\sum_{g \in K}
	(tr(\chi(\overline{g})))g^{-1}$, where $\overline{g}$ denotes the image of $g$ in $K/H,$ $\zeta_{n}$ a primitive $n^{\mathrm {th}}$ 
	root of unity in $\overline{\mathbb{F}}_{q}$ and $tr:=tr_{\mathbb{F}_{q}(\zeta_{n})/\mathbb{F}_{q}}$ denotes the trace of the field extension $\mathbb{F}_{q}(\zeta_{n})/\mathbb{F}_{q}$. Let $e_{C}(G,\,K,\,H)$ denote the sum of distinct $G$-conjugates of $\varepsilon_{C}(K, H).$ Broche and Rio \cite{Broche} proved that
	\begin{equation}\label{e20}
	e_{C}(G,K,H) \mathrm {~ is~ a ~primitive~ central~ idempotent ~of}~ \mathbb{F}_{q}G\end{equation}
	\noindent and
	\begin{equation}\label{E30} \mathbb{F}_{q}Ge_{C}(G,K,H)\cong M_{[G:K]}(\mathbb{F}_{q^{o/[E
			:K]}}), \end{equation}
		\noindent	where $E=E_{G}(K/H)$ and $o=o_{[K:H]}(q)$.\\
	
	\subsection{Main theorem}
		The following theorem gives a complete set of primitive central idempotents of finite semisimple group algebra $\mathbb{F}_{q}G$, when $G$ is a normally monomial group. \linebreak Consequently, the complete algebraic structure of $\mathbb{F}_{q}G$ is obtained.
	
		\begin{theorem}\label{T1}Let  $\mathbb{F}_{q}$ be a finite field with $q$ elements and let 
		$G$ be a finite  group. Suppose that $\operatorname{gcd}(q,\,|G|) = 1.$ Then,
		$$E:=\{e_{C}(G,\,A,\,D)\,\mid\,(A,\,D)\in \mathcal{S}(G),\, C\in \mathcal{R}(A/D)\}$$ is a complete
		set of primitive central idempotents of $\mathbb{F}_{q}G$ if, and only if, $G$ is normally monomial.
		
	\end{theorem}

In order to prove the above theorem, we need the following lemmas:

\begin{lemma}\label{l2}{\rm (\cite{BGP}, Lemma 1)} Let $\mathbb{F}_{q}$ be a finite field with $q$ elements and let 
	$G$ be a finite  group. Suppose that $\operatorname{gcd}(q,\,|G|) = 1.$ Let $\psi$ be a linear character of a normal subgroup $A$ of $G$ with kernel $D.$ 
	If ${\psi}^{G}\in \mbox{Irr}(G),$ then $e_{\mathbb{F}_{q}} ({\psi}^{G})=e_{C}(G,\,A,\,D),$ for some$~ C\in \mathcal{C}(A/D).$
	
\end{lemma}
\begin{lemma}\label{l1}{\rm(\cite{How}, Lemma 2.2)} Let $A$ be an abelian normal subgroup of maximal order in a group $G$. If $\chi$ is a faithful irreducible normally monomial character of $G$, then $\chi$ is induced from $A$. \end{lemma}
{\bf Proof.} By assumption, $\chi = \beta^{G}$ for some linear character $\beta$ of a
normal subgroup $B$ of $G$. As $B'\subseteq \operatorname{core}_{G}(\operatorname {ker} (\beta))=\operatorname{ker} (\chi)= \langle 1\rangle$, we get that $B$ is abelian and hence $\chi(1)=[G:B] \geq [G:A].$ Let $\alpha$ be an irreducible constituent of $\chi_{A}$, where $\chi_{A}$ \linebreak denotes the restriction of $\chi$ to $A$. By Frobenius reciprocity (\cite{Mus1}, Corollary 4.2.2), $\chi$ is a constituent of $\alpha^{G}$, so that $\chi(1) \leq \alpha^{G}(1) = [G : A]$. Thus, $\chi(1) = [G : A] = \alpha^{G}(1)$ and $\chi = \alpha^{G}$.~$\Box$\\

\noindent{\bf{Proof of Theorem \ref{T1}.}} Assume first that $E$ is a complete set of primitive \linebreak central idempotents of $\mathbb{F}_{q}G.$ 
	Let $\chi \in \mbox{Irr}~(G),$ so that  $e_{\mathbb{F}_{q}} ({\chi})=e_{C}(G,\,A,\,D),$ where $(A,\,D)\in \mathcal{S}(G)$ and $C\in \mathcal{R}(A/D)$. But then, it follows from Lemma \ref{l2} that
	$e_{C}(G,\,A,\,D)= e_{\mathbb{F}_{q}} ({\psi}^{G}),$ for some linear character 
	$\psi$ of $A$ with kernel $D.$ Hence $e_{\mathbb{F}_{q}} ({\chi})=e_{\mathbb{F}_{q}} ({\psi}^{G}),$ which implies that $\chi = \sigma \circ {\psi}^{G}$, where 
	$\sigma \in \operatorname{Gal}(\mathbb{F}_{q}({\psi}^{G})/\mathbb{F}_{q})$, i.e., $\chi$ is induced from a linear character $\sigma \circ \psi$ of a normal subgroup 
	$A$ of $G.$ Thus, $G$ is normally monomial.
	
	\par Conversely, if $G$ is a normally monomial group, then 
	for $\chi \in\operatorname{Irr}(G)$,\,$\chi = {\psi}^{G},$ where $\psi$ is a linear character of some normal subgroup $K$ of $G$ 
	with kernel $H.$ Now, $\operatorname{ker}(\chi)= \operatorname{core}_{G}(H)=N $(say). Let $\overline{\chi}$ be the corresponding character of $G/N$. It follows from Lemma \ref{l1} that  \begin{equation}\label{e0} \overline{\chi} = \overline{\eta}^{G/N}\end{equation}  for some linear character $\overline{\eta}$ of $A_{N}/N $, where $A_{N}/N$ is an abelian normal \linebreak subgroup of maximal order in $G/N$. Let $\ker(\overline{\eta})= L/N$ and let $\eta : A_{N}\rightarrow \overline{\mathbb{F}}_{q}$ be given  by $\eta(g)= \overline{\eta}(gN),~g\in A_{N} $. We have, by Eq.\,(\ref{e0}),  $\chi = \eta^{G}$, and therefore, by  Lemma \ref{l2}, it follows that $e_{\mathbb{F}_{q}}(\chi)= e_{\mathbb{F}_{q}}({\eta}^{G})=e_{C}(G,\,A_{N},\,L),~ C\in \mathcal{C}(A_{N}/L).$ 
	 Now, if $D\in \mathcal{T}_{N}$ is
	a representative of conjugacy class of $L$, then $D={L}^{x}, $ for some $x\in G$ and $e_{C}(G,\,A_{N},\,L)=
	e_{C^{x}}(G,\,A_{N},\,L^{x}),$ with $(A_{N},\,D)\in 
	 \mathcal{S}(G)$.\linebreak Consequently, $e_{\mathbb{F}_{q}}(\chi)\in E.$ Moreover, it is easy to check that all elements in $E$ are distinct and hence, $E$ is a complete set of primitive central idempotents of $\mathbb{F}_{q}G$.  $\Box$\\

	 \noindent It thus follows from Theorem \ref{T1} and Eq.\,(\ref{E30}) that if $G$ is normally monomial, then 
	\begin{equation}\label{E3}\mathbb{F}_{q}G\cong \displaystyle{\bigoplus_{(A,D)\in \mathcal{S}(G)}}M_{[G:A]}{(\mathbb{F}_{q^{o/[E
				:A]}})}^{(|\mathcal{R}(A/D)|)}, \end{equation}
	
	\noindent	where $E=E_{G}(A/D)$ and $o=o_{[A:D]}(q)$.
	\subsection{An illustration}
	 \noindent We first illustrate the results of Theorem 1 on a family of metabelian groups of order $2sp^{2}$, $p$ an odd prime and $s\geq 1$. The structure of the rational group algebra of this family has been described in \cite{BM1}. We undertake the case of finite semisimple group algebra for these groups.
	
	\begin{theorem}\label{th1}Let $G$ be a group generated by $a,\,b,\,x,\,y$ satisfying 
		$a^p=b^p=x^s=y^2=1,~[a,b]=[x,y]=1$,
		$x^{-1}ax=a^r,$ $x^{-1}bx=b^r,~y^{-1}ay=b,~y^{-1}by=a$, where $p$ is an odd prime, $r$ is a positive integer such that $\operatorname{gcd}(r,\,p)=1$ and  $s(>1)$ is the order $o_{p}(r)$ of $r$ modulo $p$. The Wedderburn decomposition of  $\mathbb{F}_{q}G$ is given as follows:\\
	
		\noindent \rm{\textbf{(i)}} If $s$ is odd, then \\
		\begin{scriptsize}
			$$\mathcal{\mathbb{F}}_qG\,\cong \,\left\{
			\begin{array}{ll}
			\displaystyle\bigoplus_{d|s}\mathcal{\mathbb{F}}_{q^{o_d}}^{(\frac{2\phi(d)}{o_d})} \bigoplus M_s\left(\mathcal{\mathbb{F}}_{q^{\frac{fg}{s}}}\right)^{(\frac{2(p-1)}{fg})} 
			\bigoplus M_{2s}\left(\mathbb{F}_{q^{\frac{fg}{s}}}\right)^{(\frac{p^2-p}{2fg})}, & ~~~~~~2\nmid f\\
			\displaystyle\bigoplus_{d|s}\mathbb{F}_{q^{o_d}}^{(\frac{2\phi(d)}{o_d})} \bigoplus 
			M_s\left(\mathbb{F}_{q^{\frac{fg}{s}}}\right)^{(\frac{2(p-1)}{fg})} \bigoplus M_{2s}\left(\mathbb{F}_{q^{\frac{fg}{2s}}}\right)^{(\frac{p-1}{fg})}\bigoplus 
			M_{2s}\left(\mathbb{F}_{q^{\frac{fg}{s}}}\right)^{(\frac{(p-1)^2}{2fg})}, & ~~~~~~2|f
			\end{array};
			\right.
			$$
		\end{scriptsize}
		\begin{flushleft}

			\rm{\textbf{(ii)}} If $s$ is even, then \\
			\begin{scriptsize}
				$$\mathbb{F}_qG\,\cong \,\left\{
				\begin{array}{ll}
				\displaystyle\bigoplus_{d|s}\mathbb{F}_{q^{o_d}}^{(\frac{2\phi(d)}{o_d})} \bigoplus M_s\left(\mathbb{F}_{q^{\frac{fg}{s}}}\right)^{(\frac{2(p-1)}{fg})} \bigoplus 
				M_{2s}\left(\mathbb{F}_{q^{\frac{fg}{s}}}\right)^{(\frac{{(p-1)}^{2}}{2fg})}\bigoplus 
				M_{s}\left(\mathbb{F}_{q^{\frac{2fg'}{s}}}\right)^{(\frac{{(p-1)}}{fg'})}, & ~~~2\nmid f\\
			\displaystyle\bigoplus_{d|s}\mathbb{F}_{q^{o_d}}^{(\frac{2\phi(d)}{o_d})} \bigoplus 
				M_s\left(\mathbb{F}_{q^{\frac{fg}{s}}}\right)^{(\frac{4(p-1)}{fg})} \bigoplus M_{2s}\left(\mathbb{F}_{q^{\frac{fg}{s}}}\right)^
				{(\frac{{(p-1)}^{2}}{2fg})}, & ~~~2|f
				\end{array},
				\right.
				$$
				
			\end{scriptsize}
		\end{flushleft}
		
		\vspace{.25 cm}
		
		\noindent where
		$o_d=o_d(q)$, $f=o_p(q)=o_{2p}(q),$ $g=\operatorname{gcd}\{d|s : \frac{s}{d}|f\}$ and \linebreak $g'=\operatorname{gcd}\{d\,|\,\frac{s}{2}\,:\,
		\frac{s}{2d}~\mbox{is odd and},~\frac{s}{2d}\,|\,f\}.$
		
	\end{theorem}

	\noindent\textbf{Proof.} The distinct normal subgroups of $G$ (\cite{BM1}, Lemma 2) are  $N_{0}:=\langle 1 \rangle,$ $N_{1}:=\langle ab\rangle,$~ $N_{2}:=\langle ab^{-1}\rangle,$~ $N_{3}:=\langle ab^{-1}, y\rangle,$ $G_{d}:= \langle a, b, x^{d}\rangle,$  $H_{d}:= \langle a, b, x^{d}, y\rangle,$ $d|s,$ if $s$ is odd; and, in addition, $ N_{4}:= \langle ab, x^{\frac{s}{2}}y\rangle,$  $K_{d}:= \langle a, b, x^{\frac{d}{2}}y \rangle,~$ with $d$ even and $ d|s,$ if $s$ is even. Further, for each normal subgroup $N$ of $G$, the corresponding set $\mathcal{S}_{N}$ of strong Shoda pairs of $G$ (\cite{BM1}, Lemma 3) is as follows:\\	
	
	\noindent (i) $\mathcal{S}_{N_{0}} =\{(\langle a,b \rangle, \langle  b \rangle)\} \cup \{( \langle a, b \rangle, \langle ba^{\lambda ^{i}} \rangle)\, \mid \, 1\leq i \leq \frac{p-3}{2}\}, $ where $\lambda$ is generator of the multiplicative group  of reduced residue classes modulo $p$.\vspace{.35cm}\\
	(ii) $\mathcal{S}_{N_{1}} = \begin{cases} \{(G_{s},N_{1})\},&  s ~{\rm odd} \vspace{.35cm}\\ \{( K_{s},N_{1})\},& s~{\rm even.}\end{cases}$ \vspace{.35cm}\\ 
	(iii) $\mathcal{S}_{N_{i}} =  \{(H_{s},N_{i})\}, ~ i =2,~3$. \vspace{.35cm}\\(iv)
	$\mathcal{S}_{N_{4}} =  \{( K_{s},N_{4})\}.$  \vspace{.35cm}\\
	(v) $\mathcal{S}_{G_{d}} = \begin{cases} \{( G,G_{d})\},&  d ~{\rm odd}\\  
	\emptyset  ,& d~{\rm even.}\end{cases}$ \vspace{.35cm}\\(vi)  $\mathcal{S}_{H_{d}}= \{ ( G,H_{d})\},~d|s.$ \vspace{.35cm}\\
	(vii)
	$\mathcal{S}_{K_{d}}= \{ (G,K_{d})\},~d|s$, $d$ {\rm even}. \\
	
	\noindent We tabulate the computations of parameters involved in applying Theorem \ref{T1}.\\
	
	\textbf{$s$ odd}
		\begin{small}	\[\arraycolsep=.5pt\def\arraystretch{1.4}
		\begin{tabular}{|l|c|c|c|}
			\hline
			$(A_N, D)\in \mathcal{S}(G)$ & $E_G(A_N,D)$ & $o(A_N,D)$ & $|\mathcal{R}(A_N/D)|$\\
			\hline
			$(G,G_d),~d|s$ & $G$& $o_{2d}(q)=o_d(q)$ & $\frac{\phi(d)}{o_d(q)}$\\
			\hline
			$(G,H_d), ~d|s$ & $G$ & $o_d(q)$ & $\frac{\phi(d)}{o_d(q)}$\\
			
			\hline
			$(G_s,N_1)$ & $\left\{\begin{array}{ll}<a,b,x^g,y>,& ~2|f\\<a,b,x^g>, & ~2\nmid f\end{array}\right.$ & $\left\{\begin{array}{ll}\frac{fg}{2s}, & ~2|f\\\frac{fg}{s}, &~ 2\nmid f\end{array}\right.$ &
			$\left\{\begin{array}{ll}\frac{p-1}{fg}, &~ 2|f\\\frac{p-1}{2fg}, &~ 2\nmid f\end{array}\right.$\\
			\hline
			$(H_s,N_2)$ & $<a,b,x^g,y>$ & $\frac{fg}{s}$ & $\frac{p-1}{fg}$\\
			\hline
			$(H_s,N_3)$ & $<a,b,x^g,y>$ & $\frac{fg}{s}$ & $\frac{p-1}{fg}$\\
			\hline
			$(<a,b>,<b>)$ & $<a,b,x^g>$ & $\frac{fg}{s}$ & $\frac{p-1}{fg}$\\
			\hline
			$(<a,b>,<ba^{\lambda^i}>)$ &\multirow{2}*{$<a,b,x^g>$}   & \multirow{2}*{$\frac{fg}{s}$} & \multirow{2}*{$\frac{p-1}{fg}$}\\
			$1\leq i \leq \frac{p-3}{2}$ &  & & \\
			\hline
		\end{tabular}\]
	\end{small}

	\textbf{$s$ even}
	\begin{small}	\[\arraycolsep=.5pt\def\arraystretch{1.4}
		\begin{tabular}{|l|c|c|c|}
			\hline
			$(A_N,D)$ & $E_G(A_N,D)$ & $o(A_N,D)$ & $|\mathcal{R}(A_N/D)|$\\
			\hline
			$(G,G_d),~d|s~d~\mbox{odd}$ &  $G$ &  $o_{2d}(q)=o_d(q)$ &  $\frac{\phi(d)}{o_d(q)}$\\

			\hline$(G,K_d),~d|s,~d~\mbox{even}$ & $G$ &   $o_d(q)$ &  $\frac{\phi(d)}{o_d(q)}$\\

			\hline	
			$(K_s,N_1)$ & $\left\{\begin{array}{ll}<a,b,x^{g},y>, & ~2|f\\<a,b,x^{g'}y>, &~ 2\nmid f\end{array}\right.$ & $
			\left\{\begin{array}{ll}\frac{fg}{s}, & ~2|f\\
			\frac{2fg'}{s}, &~2\nmid f\end{array}\right.$ &
			$\left\{\begin{array}{ll}\frac{p-1}{fg}, &~ 2|f\\\frac{p-1}{2fg'}, &~ 2\nmid f\end{array}\right.$\\
			
			\hline
			
			$(K_s,N_4)$ & $\left\{\begin{array}{ll}<a,b,x^{g},y>, &~ 2|f\\<a,b,x^{g'}y>, & ~2\nmid f\end{array}\right.$ & $
			\left\{\begin{array}{ll}\frac{fg}{s}, & ~2|f\\
			\frac{2fg'}{s}, & ~2\nmid f\end{array}\right.$ &
			$\left\{\begin{array}{ll}\frac{p-1}{fg}, & ~2|f\\\frac{p-1}{2fg'}, & ~2\nmid f\end{array}\right.$\\
			
			\hline
			
			$(H_s,N_2,)$ & $<a,b,x^g,y>$ & $\frac{fg}{s}$ & $\frac{p-1}{fg}$\\
			\hline
			$(H_s,N_3)$ & $<a,b,x^g,y>$ & $\frac{fg}{s}$ & $\frac{p-1}{fg}$\\
			\hline
			$(<a,b>,<b>)$ & $<a,b,x^g>$ & $\frac{fg}{s}$ & $\frac{p-1}{fg}$\\
			\hline
			$(<a,b>,<ba^{\lambda^i}>)$ & \multirow{2}*{$<a,b,x^g>$ }& \multirow{2}*{$\frac{fg}{s}$} & \multirow{2}*{$\frac{p-1}{fg}$}\\
			$1\leq i \leq \frac{p-3}{2}$ &  & & \\
			\hline
		\end{tabular}\]
	
	\end{small}
	\vspace{0.5cm}

\noindent	Consequently, Theorem \ref{th1} follows using Eq.\,(\ref{E3}). $\Box$
\section {Applications on p-groups} 
	\subsection{A family of normally monomial groups of order $p^{7}$}
	We now illustrate the main result of this paper on a normally monomial group of order $p^{7}$,  $p$ a prime, $p\geq 5$, which is not metabelian.
	
	\begin{theorem}\label{t1}Let $p\geq 5$ be a prime and let $G$ be the group generated by $a, b, c, d, s, r$ with the following defining relations:
		\begin{quote}  $a^{p}=b^{p}=c^{p}=d^{p}=s^{p}=r^{p}=1,$\\
			$ [a,c]=[b,c]=[a,d]=[b,d]=1,~[a,b]=[c,d], $\\
			$s^{-1}as=a, ~s^{-1}bs=bc,~ ~s^{-1}cs=acd,~ s^{-1}ds= ad,$ \\
			$ [s,r]=[a,b],~[r,a]=[r,b]=[r,c]=[r,d]=1.$		
			
		\end{quote}	
		Then, 
	
			$$\mathbb{F}_{q}G\,\cong \,\mathbb{F}_{q}\bigoplus {\mathbb{F}_{q^{f}}}^{(\frac{p^3-1}{f})}\bigoplus M_{p}( {\mathbb{F}_{q^{f}}})^{(\frac{p^4-p}{f})}
			\bigoplus M_{p^{3}}( {\mathbb{F}_{q^{f}}})^{(\frac{p-1}{f})},$$
			 where $f=o_{p}(q).$

	\end{theorem}
	
	\noindent \textbf{Proof.} Let $G$ be as in statement of the theorem.  How \cite{How} proved that $G$ is a normally monomial group, which is not metabelian. We find $\mathcal{S}(G)$ to write the Wedderburn decomposition of $\mathbb{F}_{q}G$.
	
	\noindent First of all, we compute the set $\mathcal{N}$ of normal subgroups of $G$. For this,
	we begin by observing that if $\langle 1\rangle\neq N\in \mathcal{N}$, then \begin{equation}\label{e1}
	[a,b]\in N.
	\end{equation}Let $1\neq g=a^{\alpha}b^{\beta}c^{\gamma}d^{\delta}r^{\rho}s^{\zeta}[a,b]^{\eta}\in N $, where $0\leq \alpha,~\beta,~\gamma,~\delta,~\rho,~\zeta,~\eta < p.$ Since $r^{-1}gr=g[a,b]^{\zeta}\in N$, we obtain that $\zeta=0$, if $[a,b]\not \in N$. Therefore, \linebreak $g=a^{\alpha}b^{\beta}c^{\gamma}d^{\delta}r^{\rho}[a,b]^{\eta}$ and $b^{-1}gb=g[a,b]^{\alpha}$ now yields $\alpha=0$, if $[a,b]\not \in N$. \linebreak Proceeding similarly, we obtain that if $[a,b]\not\in N$, then $g=1$, which contradicts the assumption and proves (\ref{e1}).
	
	\noindent We next check that if $a\not\in N$, then \begin{equation}\label{e2}
	N=\langle[a,b]\rangle ~or ~\langle[a,b], ra^{i}\rangle ,~\mathrm{ where}~ 0\leq i<p.
	\end{equation} 
	In view of (\ref{e1}), we already have that either $N=\langle[a,b]\rangle$ or there exists an element $ g=a^{\alpha}b^{\beta}c^{\gamma}d^{\delta}r^{\rho}s^{\zeta}\in N\setminus\langle [a,b]\rangle $, where $0\leq \alpha,~\beta,~\gamma,~\delta,~\rho,~\zeta < p$. Observe that $d^{-1}gd=g[a,b]^{\gamma}a^{-\zeta}\in N$ implies $\zeta=0$, if $a\not\in N$. Continuing in this manner, repeated conjugacy by $s$ implies $\delta=\gamma=\beta=0$, if $a\not \in N$ and $g=a^{\alpha}r^{\rho}$, $\rho \neq 0$.  Thus, (\ref{e2}) follows.
	
	\noindent Next, assume  $a\in N$, so that $\langle a, b^{-1}ab\rangle (\in \mathcal{N})\leq N$, i.e., we look for normal subgroups $N$ of $G$ containing $\langle a, b^{-1}ab\rangle$, which are in one-one correspondence with the normal subgroups $\overline{N}$ of $\overline{G}:=G/\langle a, b^{-1}ab\rangle$. Proceeding as previously, we see that if $\overline{1}\neq \overline{g}=\overline{b}^{\beta}\overline{c}^{\gamma}\overline{d}^{\delta}\overline{r}^{\rho}\overline{s}^{\zeta}\in \overline{N} $, for some $0\leq \beta,~\gamma,~\delta,~\rho,~\zeta < p,$ then either $\overline{d}\in \overline{N}$ or $\overline{g}=\overline{d}^{\delta}\overline{r}^{\rho}$, $\rho \neq 0$, and in the later case, \begin{equation}\label{e3}
	N=\langle a, b^{-1}ab\rangle ~or ~\langle a, b^{-1}ab, rd^{i}\rangle , ~\mathrm{ where}~ 0\leq i<p.
	\end{equation} 
	Now, if $\overline{d}\in \overline{N}$, then $\overline{N}$ is a normal subgroup of $\overline{G}$, which contains the normal subgroup $\langle \overline{d}\rangle$ of $\overline{G}$. Hence, $\overline{N}$ corresponds to a normal subgroup $\overline{\overline{N}}$ of $\overline{\overline{G}}:=\overline{G}/\langle \overline{d}\rangle$. Continuing this process, the set of normal subgroups of $G$, is obtained as 
	\begin{equation}\label{e4}
	\mathcal{N}=\{G,N_{0},N_{1},N_{2},N_{3},N_{4},N_{5},N_{6},G^{(1)}_{i},G^{(2)}_{i},G^{(3)}_{i},G^{(4)}_{i},G^{(5)}_{i},H_{ij},K_{ij}~| ~0\leq i,j < p\}, 
	\end{equation} 
	where  \begin{quote}$N_{0}:=\langle 1 \rangle; ~N_{1}:=\langle a,b,c,d,r \rangle; N_{2}:=\langle a,b,c,d \rangle;~N_{3}:=\langle a,c,d \rangle;$\\$N_{4}:=\langle a,b^{-1}ab, d \rangle;~N_{5}:=\langle a,b^{-1}ab \rangle;~N_{6}:=\langle [a,b] \rangle;$\\$G^{(1)}_{i}:=\langle a,b,c,d,sr^{i} \rangle;~G^{(2)}_{i}:=\langle a,c,d,rb^{i} \rangle;~G^{(3)}_{i}:=\langle a,b^{-1}ab,d,rc^{i} \rangle;$\\$G^{(4)}_{i}:=\langle a,b^{-1}ab,rd^{i} \rangle;~G^{(5)}_{i}:=\langle [a,b],ra^{i} \rangle;$\\$H_{ij}:=\langle a,c,d,rb^{i},sb^{j} \rangle;$ and $K_{ij}:=\langle a,c,d,sr^{i}b^{j} \rangle$, for $~0\leq i,j < p$.\end{quote}
	
	Note that if $N=G,N_{1},H_{ij},G^{(1)}_{i},N_{2},K_{ij},G^{(2)}_{i} $ or $N_{3}$, then $G/N$ is abelian and hence, by Eq.\,(\ref{e13}), we have
	
	\begin{equation}\label{e5} \mathcal{S}_{N} = \begin{cases} \{(G, N)\}, & {\rm if}~ N\in\mathcal{N}_{1}:=\{G,N_{1},H_{ij},G^{(1)}_{i}\},\\ \emptyset, & {\rm if}~N\in\mathcal{N}_{2}:=\{N_{2},K_{ij},G^{(2)}_{i},N_{3}\}.
	\end{cases} \end{equation}
	Further, if $N\neq N_{0}$ is such that $G/N$ is not abelian, then $N_{1}/N$ is abelian.\linebreak Consequently,

	\begin{equation}\label{e6} A_{N}=N_{1}, ~{\rm if} ~N\in\mathcal{N}_{3}:=\{G^{(3)}_{i},N_{4},G^{(4)}_{i},N_{5},G^{(5)}_{i},N_{6}\}.\end{equation}
	Also, it is easy to observe that \begin{equation}\label{e7}A_{N}=G^{(3)}_{0}, ~{\rm if }~N\in\mathcal{N}_{4}:=\{N_{0}\}\end{equation} serves the purpose. \\
	
	Since $\mathcal{N}=\mathcal{N}_{1} \,\cup\,\mathcal{N}_{2} \,\cup\,\mathcal{N}_{3} \,\cup\,\mathcal{N}_{4}$, in view of (\ref{e5})-(\ref{e7}), we need to find $\mathcal{T}_{N}$,  when $
	N\in \mathcal{N}_{3} \,\cup\,\mathcal{N}_{4}$. We first take up the case when $N\in \mathcal{N}_{4}$, i.e., $N=N_{0}=\langle 1\rangle$. In this case, since $A_{N_{0}}=G^{(3)}_{0}$ is an abelian $p$-group, the set $\mathcal{C}_{0}$ of subgroups of  $A_{N_{0}}$, which give cyclic quotient can be computed using the algorithm given in Section~5 of \cite{Basm} and is as follows:
	$$\mathcal{C}_{0}=\{H^{(0)},H^{(1)},H_{\beta},H_{\alpha\beta},H_{\alpha\beta\delta}~|~0\leq \alpha,\beta,\delta<p\},$$
	where 
	$H^{(0)}:=\langle a, d, r , b^{-1}ab\rangle,~H^{(1)}:=\langle a, d, r  \rangle,~H_{\beta}:=\langle d, r , b^{-1}aba^{\beta}\rangle,$ \linebreak
	 $~H_{\alpha\beta}:=\langle ad^{\alpha}, r , b^{-1}abd^{\beta}\rangle,~H_{\alpha\beta\delta}:=\langle ar^{\alpha}, b^{-1}abr^{\beta}, dr^{\delta} \rangle,$ for $0\leq \alpha,\,\beta,\,\delta<p$.  In view of (\ref{e1}), we note that for $H\in \mathcal{C}_{0}$, $\operatorname{core}_{G}(H)=N_{0}=\langle 1 \rangle$, if, and only if, $[a,b]\not\in H $. This yields,
	\begin{equation}\label{e8} \mathcal{D}_{N_{0}}=\{H^{(1)}\} \,\cup\,\{H_{\beta}~|~\beta\neq -1\}\cup\{H_{\alpha\beta},~H_{\alpha\beta\delta}~|~\alpha\neq\beta,~0\leq \alpha,\beta,\delta<p\}.\end{equation}
	It turns out that all subgroups of $\mathcal{D}_{N_{0}}$ are $G$-conjugates and hence 
	\begin{equation}\label{e9}
	\mathcal{T}_{N_{0}}=\{H^{(1)}\}.
	\end{equation}
	For ease of the reader, we list the conjugating element $g\in G$ such that $H^{g}=H^{(1)}$, for every $H\in \mathcal{D}_{N_{0}}$.

	\begin{itemize}
		\item $H_{\beta}^{x}=H^{(1)}$, with $x=b^{\nu}$, where $\nu(\beta+1)\equiv -1
		\,(\mathrm{mod}~ p)$, $\beta\neq -1$.
		\item $H_{\alpha\beta}^{y}=H^{(1)}$, with $y=(b^{\alpha}c)^{\lambda}$, where $\lambda(\beta-\alpha)\equiv 1\,({\mathrm{mod}~ p})$, $\alpha\neq \beta$.
		\item $H_{\alpha\beta\delta}^{z}=H^{(1)}$, with $z=(b^{\alpha}c^{-\delta})^{\lambda}s^{\lambda}$, where $\lambda(\beta-\alpha)\equiv 1\,({\mathrm{mod}~ p})$, $\alpha\neq \beta$.	
		
	\end{itemize}
	
	Finally, we proceed to find $\mathcal{S}_{N}$, when $N\in \mathcal{N}_{3}$. For this, we find the set $\mathcal{C}_{1}$ of subgroups of $A:=N_{1}=\langle a, b , c, d , r\rangle$ ( see (\ref{e6})) which give cyclic quotient, so that 
	\begin{equation}\label{e10}
	\mathcal{D}_{N}:=\{D\in \mathcal{C}_{1}~|~\operatorname{core}_{G}(D)=N\}.
	\end{equation}
	Let $D\in \mathcal{C}_{1}$. Then, $A/D$ being cyclic ensures $[a,b]\in D$. Moreover, any element $\overline{g}=\overline{a}^{\alpha}\overline{b}^{\beta}\overline{c}^{\gamma}\overline{d}^{\delta}\overline{r}^{\rho}\in \overline{A}:=A/D $, where $0\leq \alpha,~\beta,~\gamma,~\delta,~\rho< p,$ is such that $\overline{g}^{p}=\overline{1}$, so that $\overline{A}=\langle \overline{1}\rangle$ or a cyclic group of order $p$. If $a\not\in D$, then $\overline{a}\in \overline{A}$ yields $ \overline{A}=\langle \overline{a}\rangle.$ Also, $\overline{b}\in \overline{A}$ implies $ba^{\beta}\in D$, for some $0\leq \beta<p$. Similarly, $ca^{\gamma},~ra^{\rho},~da^{\delta}\in D$, for some $0\leq \gamma,~\rho,~\delta<p$. This, along with (\ref{e1}) and order considerations yield $D=D_{\rho\delta\gamma\beta}:=\langle [a,b],ra^{\rho},da^{\delta},ca^{\gamma},ba^{\beta}\rangle$,  for some $0 \leq \rho,\delta,\gamma,\beta<p$. Assuming next, $a\in D$ and $b\not\in D$ gives $D=D_{\rho\delta\gamma}:=\langle a,[a,b],rb^{\rho},db^{\
		\delta}
	,cb^{\gamma}\rangle$,~$0 \leq \rho,\delta,\gamma<p$. Iterating in the above manner, we obtain
	\begin{equation}\label{e11}
	\mathcal{C}_{1}:=\{D^{(0)},D^{(1)},D_{\rho},D_{\rho\delta},D_{\rho\delta\gamma},D_{\rho\delta\gamma\beta}| ~0\leq \rho,\delta,\gamma,\beta<p\},
	\end{equation}
	where $D^{(0)}:=A,~D^{(1)}:=\langle a,b,c,d \rangle$,  $D_{\rho}:=\langle a,b,c,rd^{\rho}\rangle$, $D_{\rho\delta}:=\langle a,b,dc^{\delta},rc^{\rho}\rangle,$ and  $D_{\rho\delta\gamma}$, $D_{\rho\delta\gamma\beta}$ are as defined above, $0 \leq \rho,\delta,\gamma,\beta<p$. Since $\operatorname{core}_{G}(D)$ is the largest element of $\mathcal{N}$,  contained in $D$, in view of (\ref{e4}), (\ref{e10}) \& (\ref{e11}),  we obtain:
	
	\begin{itemize}
		\item $\mathcal{D}_{N}=\emptyset$, if $N=N_{4},N_{5}$ or $N_{6}$.
		\item $\mathcal{D}_{N}=\{ D_{i0}\}\cup\{ D_{\rho0\gamma}~|~ \rho\equiv-\gamma i\,(\mathrm {mod} ~p), ~0< \gamma<p\}$,  if~ $N=G^{(3)}_{i}$,~ $0\leq i < p$.
		
		\item $\mathcal{D}_{N}=\{ D_{i}\}\cup\{D_{\rho\delta},~ D_{\rho\delta \gamma}~| ~\rho\equiv-\delta i\,(\mathrm {mod}~ p),~0<\delta<p,~0\leq \gamma<p\}$,  if~ $N=G^{(4)}_{i}$,~  $0\leq i < p$.
		\item $\mathcal{D}_{N}=\{D_{i\delta\gamma\beta}~|~0\leq \delta,\gamma,\beta<p\}$, if $N=G^{(5)}_{i}$,  $0\leq i < p$.
	\end{itemize}
	We now omit the details for the remaining part and enlist the set  $\mathcal{T}_{N}, $ when$ \linebreak N=G^{(3)}_{i}$, $G^{(4)}_{i}$,  $G^{(5)}_{i}$, $0\leq i < p$.
	
	\begin{itemize}
		\item $\mathcal{T}_{N}=\{ D_{i0}\}$,  if~ $N=G^{(3)}_{i}$,~  $0\leq i < p$.
		
		\item $\mathcal{T}_{N}=\{ D_{i}\}\cup\{ D_{i\delta0}~| ~0<\delta<p\}$,  if~ $N=G^{(4)}_{i}$,~ $0\leq i < p$.
		
		\item $\mathcal{T}_{N}=\{D_{i0\gamma\beta},~0\leq \gamma, \beta<p\}$, if $N=G^{(5)}_{i}$, $0\leq i < p$.
	\end{itemize}
	In view of (\ref{E1}) and Eq.\,(\ref{E3}), the information gathered above yields the complete set of strong Shoda pairs of $G$, and their respective simple components, as tabulated below:
	
	\begin{scriptsize}
		\[\arraycolsep=6.5pt\def\arraystretch{1.41}
		\begin{tabular} {|c|c|l|c|c|}
		\hline
		\textbf{Normal Subgroup, $N$} &
		\textbf{$A_{N}$} & \textbf{$\mathcal{S}_{N}$}
		& $|\mathcal{R}(A_{N}/D)|$& $\displaystyle{\bigoplus_{(A,D)\in \mathcal{S}_{N}}}\mathbb{F}_{q}Ge_{\mathcal{C}}(G,A_{N},D)$\\ \hline
		$G$                & $G$         & $\{(G,~G)\}$    &1 &$\mathbb{F}_{q}$\\ \hline
		
		$N_{1}=\langle a,~b,~c,~d,~r\rangle$             & $G$                & $\{(G,~N_{1})\}$                                     &$\frac{p-1}{f}$ &$\mathbb{F}_{q^{f}}$ \\ \hline
		
		$H_{i,j}=\langle a, c, d, rb^{i}, sb^{j}\rangle$ &\multirow{2}*{$G$}  & \multirow{2}*{$\{(G,~H_{ij})\}$}       &\multirow{2}*{$\frac{p-1}{f}$} &$\mathbb{F}_{q^{f}}$\\
		$0\leq i,~j < p$                              &                    &                                                                                 &&      $0\leq i,~j < p $             \\ \hline
		$G_{i}^{(1)}=\langle a,b,c,d,sr^{i} \rangle $    &\multirow{2}*{$G$}  &  \multirow{2}*{$\{(G,~G_{i}^{(1)})\}$}      &\multirow{2}*{$\frac{p-1}{f}$} &$\mathbb{F}_{q^{f}}$  \\                                                        $0\leq i < p$                                 &                    &                  &&      $0\leq i < p$                    \\  \hline
		$N_{2}= \langle a,c,d,b\rangle$                  & $G$                &    $\emptyset$                                  & -& -\\ \hline
		$K_{i,j}=\langle a, c, d, sr^{i}b^{j}\rangle$,   & \multirow{2}*{$G$} &\multirow{2}*{$\emptyset$}                       &\multirow{2}*{-} &\multirow{2}*{-}   \\
		$0\leq i,~j < p$                              &                    &                                                    & &  \\ \hline
		$G_{i}^{(2)}= \langle a, c, d, rb^{i} \rangle$ , & \multirow{2}*{$G$} & \multirow{2}*{$\emptyset$}     &\multirow{2}*{-} &\multirow{2}*{-}                  \\
		$0\leq i < p$                                 &                    &                                                  & &  \\ \hline
		$N_{3}=\langle a,c,d\rangle $                    & $G$                &    $\emptyset$                                        &- &- \\ \hline
		$G_{i}^{(3)}= \langle a, d, b^{-1}ab, rc^{i} \rangle$,& \multirow{2}*{$N_{1}$}  &\multirow{2}*{$\{(N_{1},~\langle a, b, d,  rc^{i} \rangle)\}$} &\multirow{2}*{$\frac{p-1}{f}$} &$M_{p}(\mathbb{F}_{q^{f}})$\\
		$0\leq i < p$                                 &                    &                                              & &$0\leq i < p$   \\ \hline
		$N_{4}=\langle a,b^{-1}ab,~d\rangle$             & $N_{1}$            &    $\emptyset$                           &- &- \\ \hline
		$G_{i}^{(4)}= \langle a,  b^{-1}ab, rd^{i} \rangle$,& \multirow{2}*{$N_{1}$}  &$\{(N_{1},~\langle a,b,c,rd^{i}\rangle)\}~\cup$ &\multirow{2}*{$\frac{p-1}{f}$} &$M_{p}(\mathbb{F}_{q^{f}})^{(p)}$ \\
		$0\leq i < p$                                 &                    &            $\{(N_{1},~\langle a,b^{-1}ab,c,db^{\delta},rd^{i}\rangle)~ | ~1 \leq \delta < p  \}$                                   & &$0\leq i < p$      \\ \hline
		
		$N_{5}=\langle a, b^{-1}ab\rangle$               & $N_{1}$            &  $\emptyset$                                                       & -&- \\ \hline
		$G_{i}^{(5)}= \langle [a,b], ra^{i} \rangle$,& \multirow{2}*{$N_{1}$}  &\multirow{2}*{$\{(N_{1},~\langle [a,b],d,ca^{\gamma},ba^{\beta},ra^{i}\rangle)~ | ~0 \leq \beta,\gamma  < p  \}$}  &\multirow{2}*{$\frac{p-1}{f}$} &$M_{p}(\mathbb{F}_{q^{f}})^{(p^{2})}$\\
		$0\leq i < p$                                 &                    &                                              & &$0\leq i < p$     \\ \hline
		
		$N_{6}=\langle [a,b]\rangle$            & $N_{1}$            &  $\emptyset$                                              &-&-                \\\hline
		$N_{0}=\langle 1 \rangle $                       & $G_{0}^{(3)}$
		&$ \{(\langle a,~d,~r,~b^{-1}ab\rangle,~\langle a,~d,~r\rangle)\}$ & $\frac{p-1}{f}$ &$M_{p^{3}}(\mathbb{F}_{q^{f}})$ \\ \hline
		\end{tabular}\]\end{scriptsize}
	 \vspace{1cm}
	\noindent Hence, the desired result is obtained.~$\Box$\\
	
	\noindent \textbf{Remarks.} Let $G$ be as in statement of Theorem \ref{t1}. Then, by the irredundant set of strong Shoda pairs of $G$, obtained in Theorem \ref{t1} and (\cite{BM1}, Corollary 1), it follows that 
	$$\mathbb{Q}G\cong \mathbb{Q}\bigoplus\mathbb{Q}(\zeta_{p})^{(1+p+p^{2})}\bigoplus M_{p}(\mathbb{Q}(\zeta_{p}))^{(p(1+p+p^{2}))}\bigoplus M_{p^{3}}(\mathbb{Q}(\zeta_{p})).$$
	
	\noindent A complete irredundant set of strong Shoda pairs of $G$ for the case when $p=5$ has been computed in \cite{BM1}, using \texttt{GAP} \cite{GAP}, to find the Wedderburn Decomposition of the rational group algebra $\mathbb{Q}G$. However, it may be pointed out that \texttt{GAP} package \texttt{Wedderga} \cite{wedd} fails to compute the Wedderburn decomposition of $\mathbb{Q}G$ or $\mathbb{F}_{q}G$, for any $p\geq 5$. In fact, though theoretically, \texttt{Wedderga} \cite{wedd} could handle the calculation of the Wedderburn decomposition of group algebras of groups of arbitrary size but in practice, if the order of the group is greater than $5000$ then the program may crash. \texttt{Wedderga} features  the function \texttt{StrongShodaPairs(G);} that determines a complete irredundant set of strong Shoda pairs of $G$ and
	the function \texttt{PrimitiveCentralIdempotentsByStrongSP(FG);} that computes the set of primitive central idempotents of semisimple group algebra $\mathbb{F}G$. For the case when $\mathbb{F}=\mathbb{Q}$, these functions are based on the search algorithms provided by Olivieri and del R{\'{\i}}o \cite{OliA}. Based on the work in \cite{BM1}, alternative and more efficient algorithms have been given in \cite{BM3}.  Analogously, in view of Theorem \ref{T1},  improved algorithms can be written and implemented in \texttt{GAP}, for the case when $\mathbb{F}$ is a finite~field.

	\subsection{Groups of order $p^{n}$,  $p$ prime, $n<5$}
	The groups of order $p^{n}$, $p$ prime, $n<5$ are metabelian and hence normally\linebreak monomial. Therefore, by applying Theorem \ref{T1}, we now give the explicit \linebreak structure of $\mathbb{F}_{q}G$, and its automorphism group when $G$ is a group of order $p^{n}$, $p$ prime, $n<5$. 
	The Wedderburn decomposition of $\mathbb{F}_{q}G$, 
	when $G$ is abelian is well known \cite{Broche} and the automorphism group in this case can be computed as in (\cite{BGP11},~Theorem 4). We therefore restrict to the case when $G$ is non-abelian.
	\subsubsection{Groups of order $p^{3}$}
	\textbf{\underline{$p=2$}}\\	
	\par If $p=2$, then $G$ is either $Q_{8},$ the quaternion group of order $8$ or $D_{8},$ the dihedral group of order $8$. In this case, 
	the set of primitive central idempotents and the Wedderburn decomposition of $\mathbb{F}_{q}G$ can be read from (\cite{BGP}, Examples 4.3 \& 4.4).\\
	
	\noindent In fact,
	$$\mathbb{F}_{q}Q_{8}\cong \mathbb{F}_{q}D_{8}\cong {\mathbb{F}_{q}}^{(4)} \bigoplus M_{2}(\mathbb{F}_{q}), $$
	and it follows from (\cite{BGP11}, Theorem 5) that
	$$\operatorname{Aut}(\mathbb{F}_{q}Q_{8})\cong \operatorname{Aut}(\mathbb{F}_{q}D_{8})\cong S_{4}\bigoplus SL_{2}(\mathbb{F}_{q}).$$

	\noindent \textbf{\underline{$p\neq 2$}}\\
	
	If $p \neq 2$, then $G$ is isomorphic to either $ \langle a, b \,|\, a^{p^{2}}=b^{p}=1,~ b^{-1}ab= a^{1+p}\rangle$ {\rm or} $ \langle a,b,c\, |\, a^{p}=b^{p}=c^{p}=1, ab=bac, ac=ca, bc=cb \rangle$ (\cite{Hall}, \S 4.4). The strong Shoda pairs of each of these groups, found in (\cite{BM1}, Theorems 3 \& 4) along with Theorem~\ref{T1}, Eq.\,(\ref{E3}) and  (\cite{BGP11}, Theorem 5), yield the following:

	\begin{theorem}
		Let $G$ be a non-abelian group of order $p^3$, where $p$ is an odd prime. Then, 
		$$\mathbb{F}_{q}G\cong \mathbb{F}_{q}\bigoplus {\mathbb{F}_{{q}^{f}}}^{((1+p)e)}\bigoplus M_{p}{(\mathbb{F}_{{q}^{f}})}^{(e)},$$
		and 
		$$\operatorname{Aut}(\mathbb{F}_{q}G)\cong \begin{cases} 
		({\mathbb{Z}_{f}}^{((1+p)e)}\rtimes S_{(1+p)e})\bigoplus{((SL_{p}(\mathbb{F}_{{q}^{f}})\rtimes\mathbb{Z}_{f})}^{(e)}\rtimes S_{e}), 
		& f\neq 1\\
		S_{(2+p)e} \bigoplus( SL_{p}{(\mathbb{F}_{q})}^{(e)} \rtimes S_{e}), & f=1 \\
		\end{cases},$$\\
		where $f=o_{p}(q)$ and $ e=\frac{p-1}{f}$.

	\end{theorem}

	\subsubsection{Groups of order $p^{4}$}
	
	\textbf{\underline{$p=2$}}\\
	
	If $p=2$, then up to isomorphism, there are $9$ non-isomorphic groups of order $2^{4}$ (\cite{burn}, \S 118) as listed below:
	
	\begin{description}
		
		\item $\mathcal{H}_{1}:=\langle a, b : a^{8} =b^{2} =1, ba=a^{5}b  \rangle$;
		\item $\mathcal{H}_{2}:=\langle a, b, c : a^{4} =b^{2}=c^{2} =1, cb=a^{2}bc, ab=ba,  ac=ca  \rangle$;
		\item $\mathcal{H}_{3}:=\langle a, b : a^{4} =b^{4} =1, ba=a^{3}b  \rangle$;
		\item $\mathcal{H}_{4}:=\langle a, b, c : a^{4} =b^{2} =c^{2}=1, ca=a^{3}c, ba=ab, cb=bc  \rangle$;
		\item $\mathcal{H}_{5}:=\langle a, b, c : a^{4} =b^{2} =c^{2}=1, ca=abc, ba=ab, cb=bc   \rangle$;
		\item $\mathcal{H}_{6}:=\langle a, b, c : a^{4} =b^{4} =c^{2}=1, ba=a^{3}b, ca=ac, cb=bc,  a^{2}=b^{2}  \rangle$;
		\item $\mathcal{H}_{7}:=\langle a, b : a^{8} =b^{2} =1, ba=a^{7}b  \rangle$;
		\item $\mathcal{H}_{8}:=\langle a, b : a^{8} =b^{2} =1, ba=a^{3}b  \rangle$;
		\item $\mathcal{H}_{9}:=\langle a, b : a^{8} =b^{4} =1, ba=a^{7}b, a^{4}=b^{2}  \rangle.$
	\end{description}
	
	
	\begin{theorem}\label{c21} The Wedderburn decomposition of $\mathbb{F}_{q}\mathcal{H}_{i}$, $1\leq i \leq 9$, is as follows:
		
		\begin{description}
			\item[(i)] $\mathbb{F}_{q}\mathcal{H}_{1}\cong \begin{cases} 
			{\mathbb{F}_{q}}^{(8)} \bigoplus {M_{2}(\mathbb{F}_{q})}^{(2)} & q\equiv 1,5\,({\rm{mod}}\,8)\\ 
			{\mathbb{F}_{q}}^{(4)} \bigoplus {\mathbb{F}_{{q}^{2}}}^{(2)}\bigoplus M_{2}(\mathbb{F}_{{q}^{2}}) & q\equiv 3,7\,({\rm{mod}}\,8)\\
			\end{cases}$;
			\item [(ii)] $\mathbb{F}_{q}\mathcal{H}_{2}\cong \begin{cases} 
			{\mathbb{F}_{q}}^{(8)} \bigoplus {M_{2}(\mathbb{F}_{q})}^{(2)} & q\equiv 1\,({\rm{mod}}\,4)\\ 
			{\mathbb{F}_{q}}^{(8)} \bigoplus M_{2}(\mathbb{F}_{{q}^{2}}) & q\equiv 3\,({\rm{mod}}\,4)\\
			\end{cases}$;
			\item [(iii)] $\mathbb{F}_{q} \mathcal{H}_{3} \cong \begin{cases} 
			
			{\mathbb{F}_{q}}^{(8)} \bigoplus {M_{2}(\mathbb{F}_{q})}^{(2)} & q\equiv 1\,({\rm{mod}}\,4)\\ 
			{\mathbb{F}_{q}}^{(4)} \bigoplus {\mathbb{F}_{{q}^{2}}}^{(2)}\bigoplus{M_{2}(\mathbb{F}_{q})}^{(2)}  & q\equiv 3\,({\rm{mod}}\,4)\\
			\end{cases}$;
			
			\item [(iv)] $\mathbb{F}_{q} \mathcal{H}_{4} \cong {\mathbb{F}_{q}}^{(8)} \bigoplus {M_{2}(\mathbb{F}_{q})}^{(2)}$;
			
			\item [(v)] $\mathbb{F}_{q} \mathcal{H}_{5} \cong \begin{cases}
			{\mathbb{F}_{q}}^{(8)} \bigoplus {M_{2}(\mathbb{F}_{q})}^{(2)} & q\equiv 1\,({\rm{mod}}\,4)\\ 
			{\mathbb{F}_{q}}^{(4)} \bigoplus {\mathbb{F}_{{q}^{2}}}^{(2)}\bigoplus {M_{2}(\mathbb{F}_{q})}^{(2)} & q\equiv 3\,({\rm{mod}}\,4)\\
			\end{cases}$;
			
			\item [(vi)] $\mathbb{F}_{q}\mathcal{H}_{6}\cong {\mathbb{F}_{q}}^{(8)} \bigoplus {M_{2}(\mathbb{F}_{q})}^{(2)}$;

			\item [(vii)] $\mathbb{F}_{q} \mathcal{H}_{7} \cong \begin{cases} 
			
			{\mathbb{F}_{q}}^{(4)} \bigoplus {M_{2}(\mathbb{F}_{q})}^{(3)} & q\equiv 1,7\,({\rm{mod}}\,8)\\ 
			{\mathbb{F}_{q}}^{(4)}\bigoplus M_{2}(\mathbb{F}_{{q}^{2}}) \bigoplus M_{2}(\mathbb{F}_{q}) & q\equiv 3,5\,({\rm{mod}}\,8)\\
			\end{cases}$;

			\item [(viii)] $\mathbb{F}_{q} \mathcal{H}_{8} \cong \begin{cases} 
			{\mathbb{F}_{q}}^{(4)} \bigoplus {M_{2}(\mathbb{F}_{q})}^{(3)} & q\equiv 1,3\,({\rm{mod}}\,8)\\ 
			{\mathbb{F}_{q}}^{(4)}\bigoplus M_{2}(\mathbb{F}_{{q}^{2}}) \bigoplus M_{2}(\mathbb{F}_{q}) & q\equiv 5,7\,({\rm{mod}}\,8)\\
			\end{cases}$;

			\item [(ix)] $\mathbb{F}_{q} \mathcal{H}_{9} \cong \begin{cases} 
			{\mathbb{F}_{q}}^{(4)} \bigoplus {M_{2}(\mathbb{F}_{q})}^{(3)} & q\equiv 1,7\,({\rm{mod}}\,8)\\ 
			{\mathbb{F}_{q}}^{(4)}\bigoplus M_{2}(\mathbb{F}_{{q}^{2}}) \bigoplus M_{2}(\mathbb{F}_{q}) & q\equiv 3,5\,({\rm{mod}}\,8)\\
			\end{cases}$.
		\end{description}

	\end{theorem}

	\noindent{\bf Proof.}
	(i) Define $N_{0}:= \langle 1 \rangle$, $N_{1}:= \langle a^{4} \rangle$, $N_{2}:= \langle a^{2} \rangle$, $N_{3}:= \langle a \rangle$, $H_{i}:= \langle a^{4}, a^{2i}b \rangle$, $K_{j}:= \langle a^{2}, a^{j}b \rangle$ where $0 \leq i,j \leq 1$. Observe that these subgroups are normal in $\mathcal{H}_{1}$. Using Eq.\,(\ref{e13}), we have $\mathcal{S}_{N_{1}}= \mathcal{S}_{N_{2}}= \phi$, $\mathcal{S}_{N_{3}}=\{(\mathcal{H}_{1},N_{3} )\}$, $\mathcal{S}_{H_{i}}=\{(\mathcal{H}_{1},H_{i})\}$, $\mathcal{S}_{K_{j}}=\{(\mathcal{H}_{1},K_{j})\}$, $0 \leq i,j \leq 1$. In order to find $\mathcal{S}_{N_{0}}$, we see that $\langle a \rangle$ is a maximal abelian subgroup of $\mathcal{H}_{1}$. Further, the only subgroup $D$ of $\langle a \rangle$ satisfying $\operatorname{core}_{G}(D)=\langle 1\rangle$ is $D= \langle1 \rangle$. This gives $\mathcal{S}_{N_{0}}=\{(\langle a \rangle,~\langle 1 \rangle)\}.$
	Define $$\mathcal{N}_{1}= \{\langle 1 \rangle,~ \langle a^{4} \rangle,~\langle a^{2} \rangle,~ \langle a \rangle,~\langle a,b\rangle\} ~\cup~ \{\langle a^{4}, a^{2i}b \rangle,~ \langle a^{2}, a^{j}b \rangle~|~0 \leq i,j \leq 1\}.$$ 
	
\noindent	It follows from Eqs.\,(\ref{e20})-(\ref{E30}) that $\displaystyle\bigoplus_{N\in\mathcal{N}_{1}}\bigoplus_{(A_{N},D)\in \mathcal{S}_{N}}\bigoplus_{C\in \mathcal{R}(A_{N}/D)}\mathbb{F}_{q}Ge_{C}(G,A_{N},D)$, is a direct summand of $\mathbb{F}_{q}G$ and has same $\mathbb{F}_{q}$-dimension as $\mathbb{F}_{q}G$. This yields that if $\mathcal{N}$ is the set of all normal subgroups of $\mathcal{H}_{1}$, then
	$\mathcal{S}_{N}= \phi$, if $N\not \in \mathcal{N}_{1}$, i.e., $\mathcal{S}{(\mathcal{H}_{1})}=\bigcup_{N\in \mathcal{N}_{1}}\mathcal{S}_{N}$ and
	
	$$\mathbb{F}_{q}\mathcal{H}_{1}\cong \begin{cases} 
	{\mathbb{F}_{q}}^{(8)} \bigoplus {M_{2}(\mathbb{F}_{q})}^{(2)} & q\equiv 1,5\,({\rm{mod}}\,8)\\ 
	{\mathbb{F}_{q}}^{(4)} \bigoplus {\mathbb{F}_{{q}^{2}}}^{(2)}\bigoplus M_{2}(\mathbb{F}_{{q}^{2}}) & q\equiv 3,7\,({\rm{mod}}\,8)\\
	\end{cases}.$$ 
	\noindent (ii)-(ix) For $2 \leq i \leq 9$, consider the following set $\mathcal{N}_{i}$ of normal subgroups of $\mathcal{H}_{i}$:
	\begin{small}	
		$$\begin{array}{cl}
		\mathcal{N}_{2} ~=  &\{\langle 1 \rangle,~\langle a^{2} \rangle, ~\langle a^{2},b \rangle,~\langle a,b \rangle,~\langle a,b,c \rangle \}~\cup\\&\{\langle a^{2}, b^{i}c \rangle,~\langle a, b^{i}c \rangle,~\langle ab^{i}c^{j} \rangle,~\langle a^{2}, a^{i}b, a^{j}c\rangle ~|~ 0\leq i,j \leq 1 \};\vspace{0.3cm}\\
		
		\mathcal{N}_{3} ~=  &\{\langle 1 \rangle,~\langle a^{2} \rangle,~\langle a \rangle,~\langle a,b^{2} \rangle,~\langle a,b \rangle\}~\cup\\&\{~\langle a^{2i}b^{2} \rangle,~\langle a^{2}, a^{i}b \rangle,~\langle a^{2},a^{i}b^{2} \rangle  ~|~ 0\leq i \leq 1 \};\vspace{0.3cm}\\
		\mathcal{N}_{4} ~=  &\{\langle 1 \rangle,~\langle a^{2} \rangle,~\langle a^{2},b \rangle,~\langle a,b \rangle,~\langle a,b,c \rangle \}~\cup\\&\{\langle a^{2i}b\rangle,~\langle a^{2}, b^{i}c \rangle,~\langle a, b^{i}c\rangle,~\langle a^{2}, ab^{i}c^{j} \rangle,~\langle  a^{2},a^{i}b, a^{j}c \rangle ~|~ 0\leq i,j \leq 1 \rangle\};\vspace{0.3cm}\\
		\mathcal{N}_{5} ~=  &\{\langle 1 \rangle,~\langle b \rangle,~\langle a^{2} \rangle,~\langle a^{2}b \rangle,~\langle a^{2},b \rangle, ~\langle a, b \rangle,~\langle b,ac\rangle,~\langle a^{2},b,c \rangle,~\langle a,b,c \rangle\}~\cup\\&\{\langle b,a^{2i}c \rangle ~|~ 0\leq i \leq 1 \rangle\};\vspace{0.3cm}\\
		\mathcal{N}_{6} ~=  &\{\langle 1 \rangle,~\langle a^{2} \rangle,~\langle a^{2},c \rangle,~\langle a,b \rangle,~\langle a,b,c \rangle \}~\cup\\&\{\langle a^{i}c\rangle,~\langle bc^{i}\rangle,~\langle ab^{i}c^{j}\rangle,~\langle a, b^{i}c \rangle,~\langle  a^{2},a^{i}b, a^{j}c \rangle ~|~ 0\leq i,j \leq 1 \rangle\};\vspace{0.3cm}\\
		\mathcal{N}_{7} ~=  &\{\langle 1 \rangle,~\langle a^{4} \rangle,~\langle a^{2} \rangle,~\langle a \rangle,~\langle a,b \rangle\}~\cup~ \{\langle a^{2},a^{i}b \rangle~|~ 0\leq i \leq 1 \};\vspace{0.3cm}\\
		\mathcal{N}_{8} ~= &\{\langle 1 \rangle,~\langle a^{4} \rangle,~\langle a^{2} \rangle,~\langle a \rangle,~\langle a,b \rangle\}~\cup~ \{\langle a^{2},a^{i}b \rangle~|~ 0\leq i \leq 1 \};\vspace{0.3cm}\\
		\mathcal{N}_{9} ~=  &\{\langle 1 \rangle,~\langle a^{4} \rangle,~\langle a^{2} \rangle,~\langle a \rangle,~\langle a,b \rangle\}~\cup~ \{\langle a^{2},a^{i}b \rangle~|~ 0\leq i \leq 1 \}.\\
		\end{array}$$
	\end{small}

	Proceeding as in (i), we get the following complete and irredundant set of strong Shoda pairs of $\mathcal{H}_{i}$, $2 \leq i \leq 9$, which yield the desired result. \\

	$\begin{array}{lll}
	
	(ii)~ \mathcal{S}(\mathcal{H}_{2})&\hspace{-.4cm}=&\hspace{-.4cm} \{(\mathcal{H}_{2}, ~\mathcal{H}_{2}),~~(\langle a, b\rangle, ~\langle b \rangle),~~ (\mathcal{H}_{2},~\langle a, b \rangle)\}~\cup \\
	&  &\hspace{-.4cm} \{(\mathcal{H}_{2},~\langle a,b^{i}c\rangle),~~ (\mathcal{H}_{2},~\langle a^{2},a^{i}b, a^{j}c\rangle)~| ~0 \leq i,j \leq 1 \}
	;\vspace{000.25cm}\\
	(iii)~\mathcal{S}(\mathcal{H}_{3})&\hspace{-.4cm}=&\hspace{-.4cm} \{(\mathcal{H}_{3}, ~\mathcal{H}_{3}),~~(\mathcal{H}_{3},~\langle a,b^{2}\rangle),~~(\mathcal{H}_{3},~\langle a\rangle),~~(\mathcal{H}_{3},~\langle a^{2}, ab^{2}\rangle)\}~\cup \\
	& &\hspace{-.4cm} \{(\langle a,b^{2} \rangle,~\langle a^{2i}b^{2}\rangle),~~(\mathcal{H}_{3},~\langle a^{2}, a^{i}b\rangle) ~|~0\leq i\leq 1\}
	;\vspace{000.25cm}\\
	(iv)~\mathcal{S}(\mathcal{H}_{4})&\hspace{-.4cm}=&\hspace{-.4cm} \{(\mathcal{H}_{4}, ~\mathcal{H}_{4}),~~(\mathcal{H}_{4},~\langle a, b \rangle)\} ~\cup\\
	& &\hspace{-.4cm} \{(\langle a, b \rangle,~\langle a^{2i}b\rangle),~~(\mathcal{H}_{4},~\langle a,b^{i}c\rangle),~~(\mathcal{H}_{4},~\langle a^{2},a^{i}b, a^{j}c\rangle)~|~0\leq i,j\leq 1\}
	;\vspace{000.25cm}\\
	
	(v)~\mathcal{S}(\mathcal{H}_{5})&\hspace{-.4cm}=&\hspace{-.4cm} \{(\mathcal{H}_{5}, ~\mathcal{H}_{5}),~~(\mathcal{H}_{5},~\langle a,b\rangle),~~(\mathcal{H}_{5},~\langle a^{2}, b,c \rangle),~~(\mathcal{H}_{5},\langle b,ac\rangle)
	\}~\cup \\
	& &\hspace{-.4cm}\{(\langle a, b \rangle,~\langle a \rangle),~~(\langle a, b \rangle,~\langle a^{2}b\rangle)\}
	~\cup  \{(\mathcal{H}_{5},~\langle b,a^{2i}c\rangle)~|~0 \leq i \leq 1\};\vspace{000.25cm}\\
	
	(vi)~\mathcal{S}(\mathcal{H}_{6})&\hspace{-.4cm}=&\hspace{-.4cm} \{(\mathcal{H}_{6},~\mathcal{H}_{6}),~~(\mathcal{H}_{6},~\langle a, b \rangle)\}~\cup \\
	& &\hspace{-.4cm} \{(\langle a, c \rangle,~\langle a^{i}c\rangle),~~(\mathcal{H}_{6},~\langle a,b^{i}c\rangle ),~~(\mathcal{H}_{6},~\langle a^{2},a^{i}b,a^{j}c\rangle )~|~0 \leq i,j \leq 1\};\vspace{000.25cm}\\
	
	(vii)~\mathcal{S}(\mathcal{H}_{7})&\hspace{-.4cm}=&\hspace{-.4cm} \{(\mathcal{H}_{7}, ~\mathcal{H}_{7}),~~(\mathcal{H}_{7},~\langle a\rangle)\}~\cup \\
	& & \hspace{-.4cm}\{(\langle a\rangle, ~\langle  a^{4i}\rangle),~~(\mathcal{H}_{7},~\langle  a^{2},a^{i}b \rangle)~|~0 \leq i \leq 1\}
	;\vspace{000.25cm}\\
	
	(viii)~\mathcal{S}(\mathcal{H}_{8})&\hspace{-.4cm}=&\hspace{-.4cm} \{(\mathcal{H}_{8}, ~\mathcal{H}_{8}),~~(\mathcal{H}_{8},~\langle a\rangle)\}~\cup \\
	& & \hspace{-.4cm}\{(\langle a\rangle, ~\langle  a^{4i}\rangle),~~(\mathcal{H}_{8},~\langle  a^{2},a^{i}b \rangle)~|~0 \leq i \leq 1\}
	;\vspace{000.25cm}\\
		\end{array}$
	$\begin{array}{lll}
	(ix)~\mathcal{S}(\mathcal{H}_{9})&\hspace{-.4cm}=&\hspace{-.4cm} \{(\mathcal{H}_{9}, ~\mathcal{H}_{9}),~~(\mathcal{H}_{9},~\langle a\rangle)\}~\cup \\
	& & \hspace{-.4cm}\{(\langle a\rangle, ~\langle  a^{4i}\rangle),~~(\mathcal{H}_{9},~\langle  a^{2},a^{i}b \rangle)~|~0 \leq i \leq 1\}.~\Box
	\vspace{000.25cm}\\
	\end{array}$

	\begin{cor}
		The automorphism group of $\mathbb{F}_{q}\mathcal{H}_{i}$, $1\leq i \leq 9$, is as follows:
		\begin{description}\item$ \operatorname{Aut}(\mathbb{F}_{q} \mathcal{H}_{1} )\cong \begin{cases}
			S_{8} \bigoplus (SL_{2}{(\mathbb{F}_{q})}^{(2)}\rtimes S_{2}), & 
			q\equiv 1,5\,({\rm{mod}}\,8)\\
			S_{4} \bigoplus ({\mathbb{Z}_{2}}^{(2)} \rtimes S_{2}) \bigoplus (SL_{2}(\mathbb{F}_{q^{2}}) \rtimes \mathbb{Z}_{2}), &
			q\equiv 3,7\,({\rm{mod}}\,8)\\
			\end{cases};$
			\item $ \operatorname{Aut}(\mathbb{F}_{q}\mathcal{H}_{2} )\cong \begin{cases}
			S_{8} \bigoplus (SL_{2}{(\mathbb{F}_{q})}^{(2)} \rtimes S_{2}), & q\equiv 1\,({\rm{mod}}\,4)\\
			S_{8} \bigoplus (SL_{2}(\mathbb{F}_{q^{2}}) \rtimes \mathbb{Z}_{2}), & q\equiv3\,{\rm{(mod\,4)}}\\
			\end{cases};$
			\item$ \operatorname{Aut}(\mathbb{F}_{q} \mathcal{H}_{3} )\cong \begin{cases}
			S_{8} \bigoplus (SL_{2}{(\mathbb{F}_{q})}^{(2)} \rtimes S_{2}), & q\equiv 1\,({\rm{mod}}\,4)\\
			S_{4} \bigoplus ({\mathbb{Z}_{2}}^{(2)} \rtimes S_{2}) \bigoplus 
			({SL_{2}(\mathbb{F}_{q})}^{(2)} 
			\rtimes S_{2}), &q\equiv 3\,({\rm{mod}}\,4)\\
			\end{cases};$
			\item$ \operatorname{Aut}(\mathbb{F}_{q}\mathcal{H}_{4} )\cong 
			S_{8} \bigoplus (SL_{2}{(\mathbb{F}_{q})}^{(2)} \rtimes S_{2}) 
			;$
			\item$ \operatorname{Aut}(\mathbb{F}_{q}\mathcal{H}_{5} )\cong \begin{cases}
			S_{8} \bigoplus (SL_{2}{(\mathbb{F}_{q})}^{(2)} \rtimes S_{2}), & q\equiv 1\,({\rm{mod}}\,4)\\
			S_{4} \bigoplus ({\mathbb{Z}_{2}}^{(2)} \rtimes S_{2}) 
			\bigoplus ({SL_{2}(\mathbb{F}_{q})}^{(2)} \rtimes S_{2}), &q\equiv 3\,({\rm{mod}}\,4)\\
			\end{cases};$
			\item$ \operatorname{Aut}(\mathbb{F}_{q} \mathcal{H}_{6} )\cong 
			S_{8} \bigoplus (SL_{2}{(\mathbb{F}_{q})}^{(2)}\rtimes S_{2});$
			\item$ \operatorname{Aut}(\mathbb{F}_{q} \mathcal{H}_{7} )\cong \begin{cases}
			S_{4} \bigoplus (SL_{2}{(\mathbb{F}_{q})}^{(3)} \rtimes S_{3}), & q\equiv 1,7\,({\rm{mod}}\,8)\\
			S_{4} \bigoplus (SL_{2}(\mathbb{F}_{q^{2}}) \rtimes \mathbb{Z}_{2}) \bigoplus SL_{2}(\mathbb{F}_{q}), &q\equiv 3,5\,({\rm{mod}}\,8)\\
			\end{cases};$
			\item$ \operatorname{Aut}(\mathbb{F}_{q} \mathcal{H}_{8} )\cong \begin{cases}
			S_{4} \bigoplus (SL_{2}{(\mathbb{F}_{q})}^{(3)} \rtimes S_{3}), & q\equiv 1,3\,({\rm{mod}}\,8)\\
			S_{4} \bigoplus (SL_{2}(\mathbb{F}_{q^{2}}) \rtimes \mathbb{Z}_{2}) \bigoplus SL_{2}(\mathbb{F}_{q}), &q\equiv 5,7\,({\rm{mod}}\,8)\\
			\end{cases};$
			\item$ \operatorname{Aut}(\mathbb{F}_{q} \mathcal{H}_{9} )\cong \begin{cases}
			S_{4} \bigoplus (SL_{2}{(\mathbb{F}_{q})}^{(3)} \rtimes S_{3}), & q\equiv 1,7\,({\rm{mod}}\,8)\\
			S_{4} \bigoplus (SL_{2}(\mathbb{F}_{q^{2}}) \rtimes \mathbb{Z}_{2}) \bigoplus SL_{2}(\mathbb{F}_{q}), &q\equiv 3,5\,({\rm{mod}}\,8)\\
			\end{cases}.$\end{description}

	\end{cor}
	\noindent \textbf{\underline{$p\neq 2$}}\\
	
	If $p$ is an odd prime, then, up to isomorphism, the following are non-abelian groups of order $p^{4}$ (\cite{burn}, \S 117):
	\begin{description}
		\item  $\mathcal{G}_{1}:= \langle a, b : a^{p^{3}}= b^{p}=1, ba=a^{1+p^{2}}b \rangle;$
		\item $\mathcal{G}_{2}:=\langle a, b, c: a^{p^{2}} =b^{p} =c^{p}=1, cb=a^{p}bc, ab=ba, ac=ca \rangle;$
		\item $\mathcal{G}_{3}:=\langle a, b : a^{p^{2}} =b^{p^2} =1, ba=a^{1+p}b  \rangle;$
		\item $\mathcal{G}_{4}:=\langle a, b, c: a^{p^{2}} =b^{p} =c^{p}=1, ca=a^{1+p}c, ba=ab, cb=bc \rangle;$
		\item $\mathcal{G}_{5}:=\langle a, b, c: a^{p^{2}} =b^{p} =c^{p}=1, ca=abc, ab=ba, bc=cb \rangle;$
		\item $\mathcal{G}_{6}:=\langle a, b, c: a^{p^{2}} =b^{p} =c^{p}=1, ba=a^{1+p}b, ca=abc, cb=bc \rangle;$
		\item $\mathcal{G}_{7}:=\begin{cases}\langle a, b, c: a^{p^{2}} =b^{p}=1, c^{p}=a^{p}, ab=ba^{1+p}, ac=cab^{-1}, cb=bc \rangle,~{\rm if}~ p=3, \\
		\langle a, b, c: a^{p^{2}} =b^{p} =c^{p}=1, ba=a^{1+p}b, ca=a^{1+p}bc, cb=a^{p}bc \rangle, ~{\rm if}~ p>3 ;\end{cases}$
		\item $\mathcal{G}_{8}:=\begin{cases}\langle a, b, c: a^{p^{2}} =b^{p}=1, c^{p}=a^{-p}, ab=ba^{1+p}, ac=cab^{-1}, cb=bc \rangle, ~{\rm if}~ p=3, \\
		\langle a, b, c: a^{p^{2}} =b^{p} =c^{p}=1, ba=a^{1+p}b, ca=a^{1+dp}bc, cb=a^{dp}bc \rangle,~ {\rm if}~p>3 \end{cases}$
		
		\hspace{1.5cm} $ d\not\equiv 0,1\,({\rm mod}~ p{\rm)};$
		\item $\mathcal{G}_{9}:=\langle a, b, c, d: a^{p} =b^{p} =c^{p}=d^{p}=1, dc=acd, bd=db, ad=da, bc=cb,\\~~~~~~~~~~~ac=ca,ab=ba \rangle;$
		\item $\mathcal{G}_{10}:=\begin{cases}\langle a, b, c: a^{p^{2}} =b^{p}=c^{p}=1, ab=ba, ac=cab, bc=ca^{-p}b \rangle,~ {\rm if}~ p=3,\\
		\langle a, b, c,d: a^{p} =b^{p} =c^{p}=d^{p}=1, dc=bcd, db=abd,ad=da,\\~~~~~~~~~~~~~~bc=cb,ac=ca,ab=ba \rangle, ~{\rm if}~ p >3. \end{cases}$
		
	\end{description}
	\vspace{0.4cm}
	
	For $1 \leq i \leq 10$,  a complete irredundant set  $\mathcal{S}(\mathcal{G}_{i})$ of strong Shoda pairs of $\mathcal{G}_{i}$ has been computed in (\cite{BM2}, Theorem 3) which in view of Theorem \ref{T1}, and Eq.\,(\ref{E3}) yields the following:

	\begin{theorem}\label{c22} The Wedderburn decomposition of $\mathbb{F}_{q}\mathcal{G}_{i}$, $1\leq i \leq 10$, is as follows:
		\begin{description}
			\item $\mathbb{F}_{q} \mathcal{G}_{1} \cong \mathbb{F}_{q}\bigoplus {\mathbb{F}_{{q}^{f}}}^{((1+p)e)} \bigoplus {\mathbb{F}_{{q}^{fp}}}^{(pe)}
			\bigoplus M_{p}{(\mathbb{F}_{{q}^{fp}})}^{(e)}$;
			
			\item $\mathbb{F}_{q} \mathcal{G}_{2} \cong \mathbb{F}_{q}\bigoplus {\mathbb{F}_{{q}^{f}}}^{((1+p+p^{2})e)}\bigoplus M_{p}{(\mathbb{F}_{{q}^{fp}})}^{(e)};$
			
			\item $\mathbb{F}_{q} \mathcal{G}_{3} \cong \mathbb{F}_{q}\bigoplus {\mathbb{F}_{{q}^{f}}}^{((1+p)e)} \bigoplus {\mathbb{F}_{{q}^{fp}}}^{(ep)}
			\bigoplus M_{p}{(\mathbb{F}_{{q}^{f}})}^{(pe)}$;

			\item $\mathbb{F}_{q} \mathcal{G}_{4} \cong \mathbb{F}_{q}\bigoplus {\mathbb{F}_{{q}^{f}}}^{((1+p+p^{2})e)} \bigoplus M_{p}{(\mathbb{F}_{{q}^{f}})}^{(pe)}$;
			
			\item $\mathbb{F}_{q} \mathcal{G}_{5} \cong \mathbb{F}_{q}\bigoplus {\mathbb{F}_{{q}^{f}}}^{((1+p)e)} \bigoplus {\mathbb{F}_{{q}^{fp}}}^{(pe)}	
			\bigoplus M_{p}{(\mathbb{F}_{{q}^{f}})}^{(pe)}$;

			\item $\mathbb{F}_{q} \mathcal{G}_{6} \cong \mathbb{F}_{q}\bigoplus {\mathbb{F}_{{q}^{f}}}^{((1+p)e)} \bigoplus M_{p}{(\mathbb{F}_{{q}^{f}})}^{((1+p)e)}$;
			
			\item $\mathbb{F}_{q} \mathcal{G}_{7} \cong \mathbb{F}_{q}\bigoplus {\mathbb{F}_{{q}^{f}}}^{((1+p)e)} \bigoplus M_{p}{(\mathbb{F}_{{q}^{f}})}^{(e)}	
			\bigoplus M_{p}{(\mathbb{F}_{{q}^{fp}})}^{(e)}$;
			
			\item $\mathbb{F}_{q} \mathcal{G}_{8} \cong \mathbb{F}_{q}\bigoplus {\mathbb{F}_{{q}^{f}}}^{((1+p)e)}\bigoplus M_{p}{(\mathbb{F}_{{q}^{f}})}^{(e)}
			\bigoplus M_{p}{(\mathbb{F}_{{q}^{fp}})}^{(e)}$;

			\item $\mathbb{F}_{q} \mathcal{G}_{9} \cong \mathbb{F}_{q}\bigoplus {\mathbb{F}_{{q}^{f}}}^{((1+p+p^{2})e)}\bigoplus M_{p}{(\mathbb{F}_{{q}^{f}})}^{(pe)}$;
			
			\item $\mathbb{F}_{q}\mathcal{G}_{10} \cong \left\{\begin{array}{ll}
			\mathbb{F}_{q}\bigoplus {\mathbb{F}_{{q}^{f}}}^{((1+p)e)}\bigoplus M_{p}{(\mathbb{F}_{{q}^{f}})}^{(e)}
			\bigoplus M_{p}{(\mathbb{F}_{{q}^{fp}})}^{(e)}, & \hbox{$p=3$} \\
			\mathbb{F}_{q}\bigoplus {\mathbb{F}_{{q}^{f}}}^{((1+p)e)}\bigoplus M_{p}{(\mathbb{F}_{{q}^{f}})}^{((1+p)e)}, & \hbox{$p>3$}
			\end{array},
			\right.$
			
		\end{description}
		
		\noindent where $o_{p}(q)=f$ and $ e=\frac{p-1}{f}$.
		
	\end{theorem}
	
	The automorphism group of finite semisimple group algebra of groups of order $p^{4}$, $p$ an odd prime, can now be computed similarly.\\

	\noindent\textbf{Acknowledgement}\\
	\noindent The authors are grateful to I.\,B.\,S. Passi and G.\,K.\,Bakshi for their valuable \linebreak comments and suggestions.

\bibliographystyle{amsplain}
	\bibliography{Bibliography}

\end{document}